\documentclass[10pt]{article}
\usepackage{graphicx}
\usepackage{textcomp}
\usepackage{amssymb}
\usepackage{amsmath}
\usepackage{epstopdf}
\usepackage{epsfig}
\usepackage{multicol}
\begin{document}

\begin{center}
\textbf{Practical Explicitly Invertible Approximation to 4 Decimals of Normal Cumulative Distribution Function Modifying Winitzki's Approximation of erf}
\end{center}

\noindent
\begin{center}
Alessandro Soranzo$^{\ast}$, Emanuela Epure$^{\star}$\\
$^{\ast}$
\scriptsize{Dipartimento di Matematica e Geoscienze -- University of Trieste --
Trieste -- Italy

email: soranzo@units.it}

$^{\star}$
\scriptsize{Esteco S.p.A. (Optimization Software) -- Area Science Park -- Trieste -- Italy

email: epure@esteco.com}

\end{center}

\medskip
\noindent
\textbf{Abstract.} We give a new explicitly invertible approximation of the normal cumulative distribution function:
$\Phi(x) \simeq \frac{1}{2} + \frac{1}{2} \sqrt{ 1-{e}^{-x^2\frac{17+{x}^{2}}{26.694+2x^2}}}$, $\forall x \ge 0$, with absolute error $<4.00\cdot 10^{-5}$, absolute value of the relative error $<4.53\cdot 10^{-5}$,
which, beeing designed essentially for practical use, is much simpler than a previously published formula and, though less precise, still reaches 4 decimals of precision, and has a complexity essentially comparable with that of the approximation of the normal cumulative distribution function $\Phi(x)$ immediatly derived from Winitzki's approximation of erf$(x)$, reducing about $36\%$ the absolute error and about $28\%$ the relative error with respect to that, overcoming the threshold of 4 decimals of precision.
\\

\\

\bigskip
\bigskip
\medskip

\noindent
This paper is devoted to approximate some special functions, in particular the normal cumulative distribution function. \\

\noindent
Though computers now allow to compute them with arbitrary precision, such approximations are still worth for several reasons, including \textit{to catch the soul} of the considered functions, allowing to understand at a glance their behaviour. 
Let's add that, despite technologic progress, those functions -- of wide practical use -- not always are available on pocket calculators. In practice, the ancient numerical tables are still widely used; but they give approximations only for some values and, if linearly interpolating to approximate intermediate values, both precision and simplicity of use reduce. And surely by computers you may obtain graphs of those functions, but for a mathematician the meaning content of formuls is greater.

\noindent
Furthermore, here we produce an \textit{explicitly invertible} (and, in fact, \textit{simply}) approximation, which allow to keep coherence working contemporarily with the considered function and its inverse.
\\

\newpage
\noindent
For the \cite{Wiki1} \cite{Wolfram1}
special functions

\begin{equation}
\label{normalCDF}
\Phi(x):= \int_{-\infty }^{x} \frac{1}{\sqrt {2\pi}} e^{\frac{-t^2}{2}} dt,
\end{equation}

\noindent
and the related error function 
$erf(x)=2\Phi(x\sqrt 2) - 1$ $\big( \forall x \in I\!\!R \big)$, 
there are several approximations; in particular see classical \cite{Abramowitz} \cite{Hart} \cite{Balakrishnan} and recent \cite{Dyer} \cite{SoranzoEpure} \cite{Winitzki} \cite{ZogheibHlynka}; approximations for  $x \ge 0$ are sufficient because of the symmetry formula 
$\Phi(-x)=1-\Phi (x) $ $\big ( \forall x \in I\!\!R \big)$.\\

\noindent
Restricting now our attention only to those approximations which are \textit{simply explicitly invertible} -- in the sense, explicitly invertible without requiring to solve cubic or quartic equations -- the most precise appears to be \cite{SoranzoEpure} 

\begin{equation}
\label{eq:Wini1.2735457}
\Phi(x) \cong \frac{1}{2} + \frac{1}{2} \sqrt{1-{e}^{\frac{-1.2735457x^2-0.0743968x^4}{2+0.1480931x^2+0.0002580x^4}}}
\quad
\left\{
    \begin{array}{ll}
       |\varepsilon (x)| <1.14\cdot 10^{-5}\\  
       |\varepsilon_r (x)| <1.78\cdot 10^{-5}
    \end{array}
    \forall x \ge 0
\right.
\end{equation}

\noindent
which is an improvement
preserving (despite the adding of the quartic monomial) the simple explicit invertibility (essential solving a biquadratic equation after obvious substitutions) of this approximation of $\Phi$

\begin{equation}
\label{eq:Wini735}
\Phi(x) \simeq\frac{1}{2} +\frac{1}{2}\sqrt{1-{e}^{-\frac{{x}^{2}\,\left( \frac{4}{\pi }+0.0735\,{x}^{2}\right) }{2\,\left(1+ 0.0735\,{x}^{2}\right) }}}
\quad
\left\{
    \begin{array}{ll}
       |\varepsilon (x) |<6.21\cdot 10^{-5} \\  
       |\varepsilon_r (x) |<6.30\cdot 10^{-5} 
    \end{array}
    \forall x \ge 0
\right.
\end{equation}

\noindent
immediately derived by $\Phi(x)=\frac{1}{2}+\frac{1}{2}erf \Big (\frac{x}{\sqrt 2} \Big )$ $\big (\forall x \ge 0\big )$ from this \cite{Winitzki} \textit{Winitzki's Approximation of erf} 

\begin{equation}
\label{eq:Wini147}
erf(x)\cong \sqrt{1-e^{-x^2\frac{\frac{4}{\pi}+0.147x^2}{1+0.147x^2}}}
\quad
\left\{
    \begin{array}{ll}
       |\varepsilon (x) |<1.25\cdot 10^{-4} \\  
       |\varepsilon_r (x) |<1.28\cdot 10^{-4} 
    \end{array}
    \forall x \ge 0
\right.
.
\end{equation}

\noindent
\textbf{In this note we give this new (simply) explicitly invertible approximation of the normal cumulative distribution function}
\begin{center}
\begin{equation}
\label{eq:ourPhiSimple}
\begin{tabular}{|c|}
\hline
$\Phi(x) \simeq \frac{1}{2} + \frac{1}{2} \sqrt{ 1-{e}^{-x^2\frac{17+{x}^{2}}{26.694+2x^2}}}
\quad
\left\{
    \begin{array}{ll}
       |\varepsilon (x) |<4.00\cdot 10^{-5} \\
       |\varepsilon_r (x) |<4.53\cdot 10^{-5}
    \end{array}
    \forall x \ge 0
\right.
$\\
\hline 
\end{tabular}
\end{equation}
\end{center}

\noindent
\textbf{which, beeing designed essentially for practical use,}

$\bullet$ \textbf{is much simpler than (\ref{eq:Wini1.2735457}) and, though less precise, still reaches 4 decimals of precision;}

$\bullet$ \textbf{has a complexity essentially comparable with that of (\ref{eq:Wini735}) reducing about $36\%$ the absolute error and about $28\%$ the relative error with respect to that, overcoming the threshold of 4 decimals of precision.}\\

\newpage
\noindent
Instead, the corresponding approximation of erf is not so worth, because, though reduces about $36\%$ the absolute error of (\ref{eq:Wini147}), it remains with the precision of 3 decimals, and furthermore the absolute value of the relative error
$|\varepsilon_r (x)| <1.79 \cdot 10^{-4}$ is quite greater than in (\ref{eq:Wini147}).

\noindent
See below the graphs (made by Mathematica $\textsuperscript{\textregistered}$) of the approximation, of the absolute error and of the absolute value of the relative error for $0 \le x \le 7$. For $x \ge 7$ the trivial approximation $\Phi(x) \simeq 1$ has absolute error and absolute value of the relative error highly less than $4 \cdot 10^{-5}$ and $4.53 \cdot 10^{-5}$ respectively. (Nevertheless, if interested in a formal proof of the majorization $|\varepsilon (x)| < 4 \cdot 10^{-5}$ for $x \ge 7$, you may follow \cite{SoranzoEpure}).

\begin{figure}[h]
\begin{center}$
\begin{array}{cc}
\includegraphics[width=2.2in,height=1.5in]{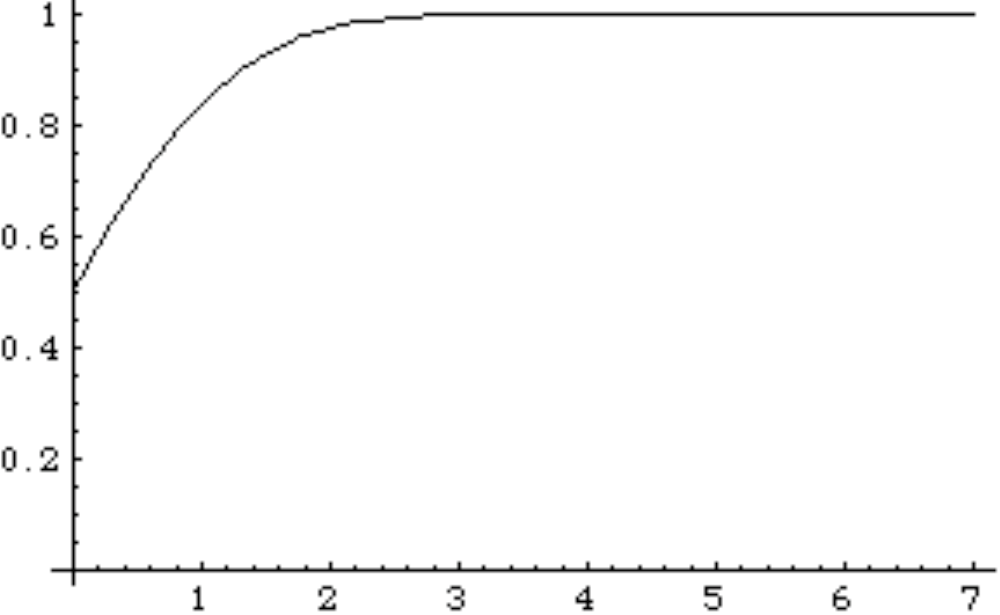} &
\includegraphics[width=2.2in,height=1.5in]{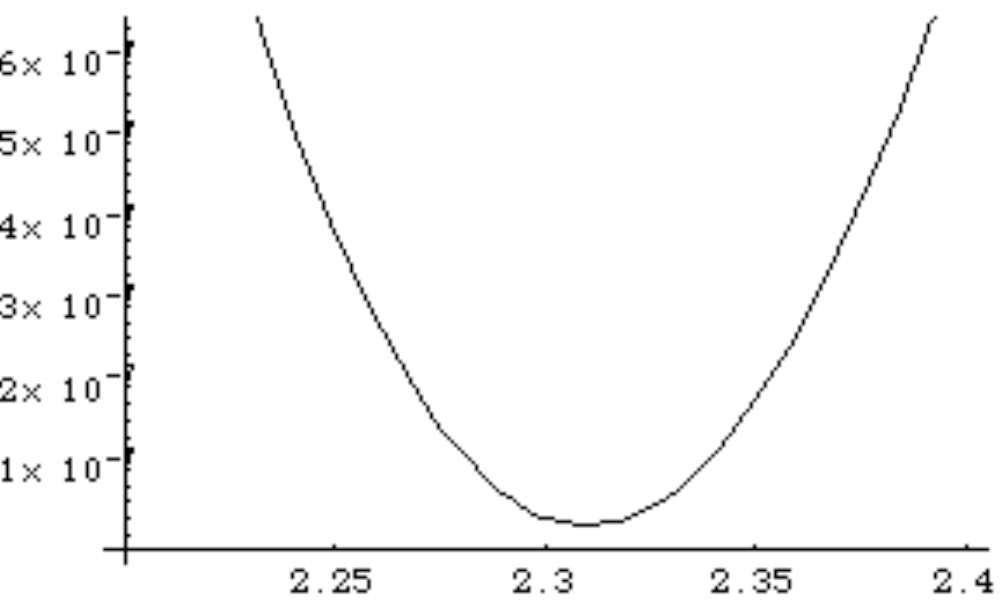} \\ 

\rm{Fig.\: 1.\: The\: new\: approximation\: (\ref{eq:ourPhiSimple})} &
\rm{Fig.\: 4.\: Second\: zoom\: of\: Fig.\: 2.}\\

\includegraphics[width=2.2in,height=1.5in]{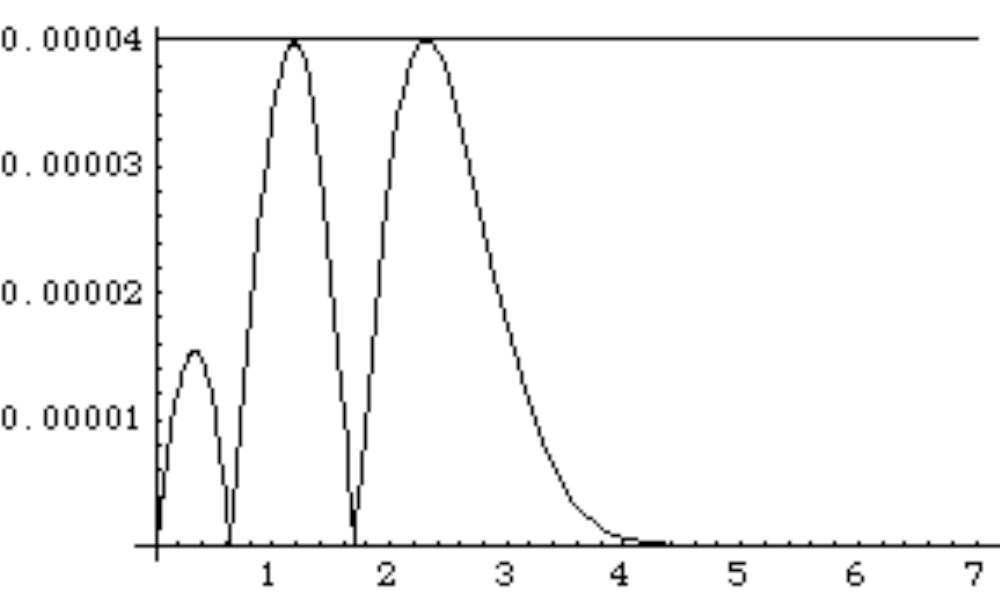} &
\includegraphics[width=2.2in,height=1.5in]{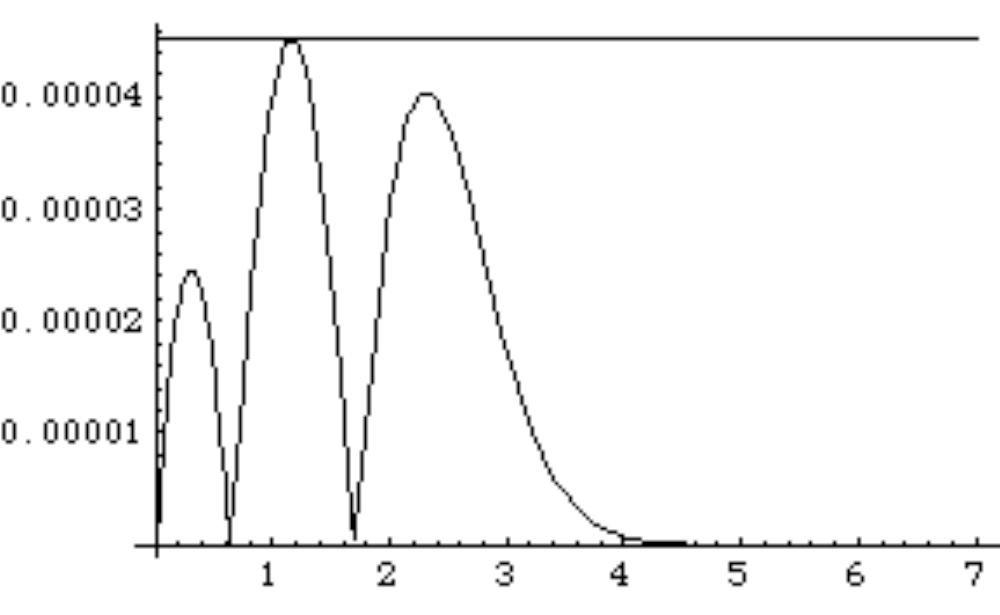} \\

\rm{Fig.\: 2.\: Absolute\: error} &
\rm{Fig.\: 5.\: Absolute\: value\: of\: relative\: error} \\

\includegraphics[width=2.2in,height=1.5in]{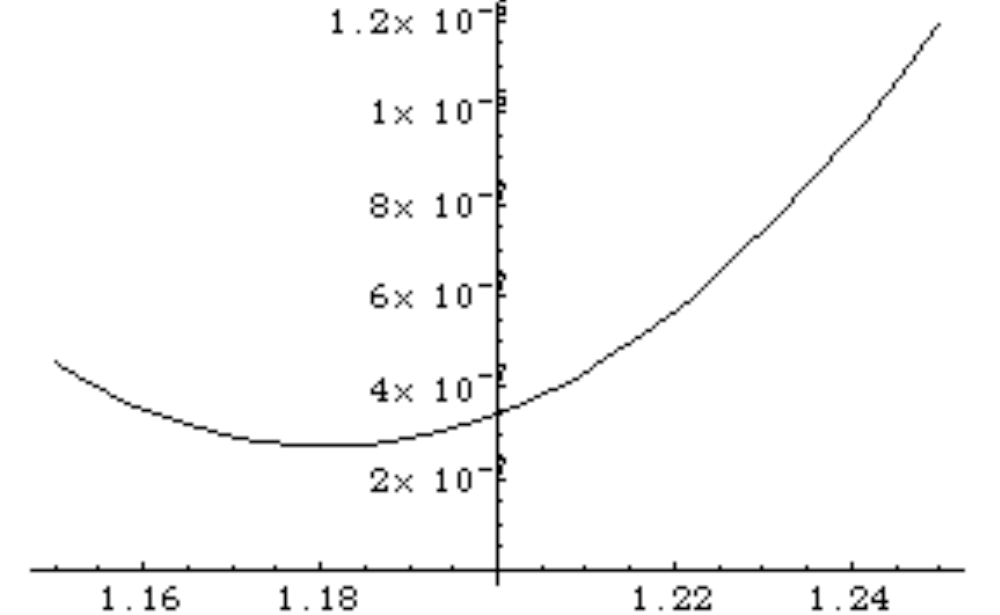} &
\includegraphics[width=2.2in,height=1.5in]{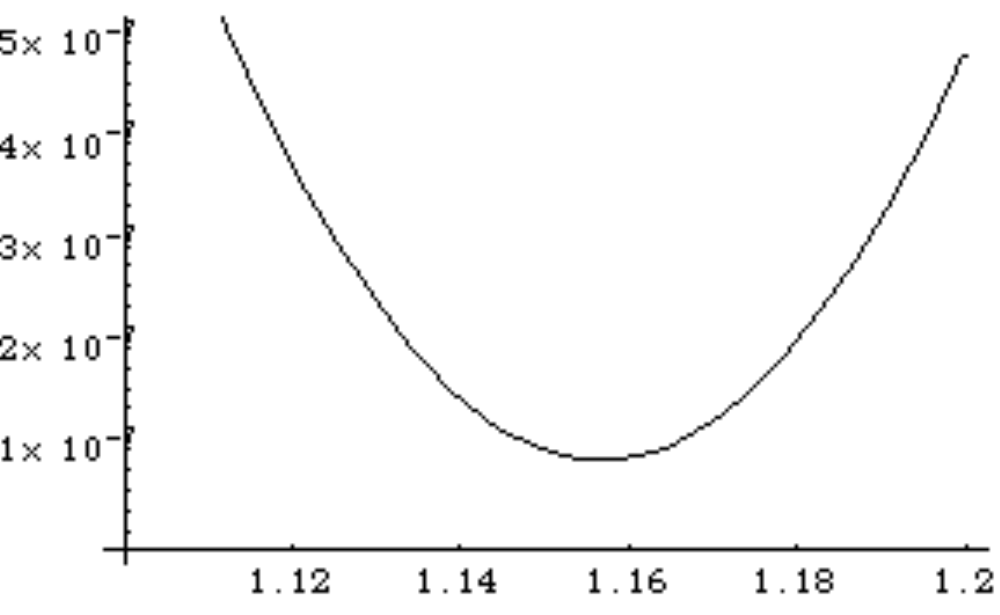} \\

\rm{Fig.\: 3.\: First\: zoom\: of\: Fig.\: 2} &
\rm{Fig.\: 6.\: Zoom\: of\: Fig.\: 5.}

\end{array}$
\end{center}
\end{figure}

\noindent
Remark. This text has been explicitly placed by the Authors in the Creative Commons Public Domain.

\end{document}